\documentclass[12pt]{amsart}
\usepackage{amssymb}

\newtheorem{Thm}{Theorem}[section]

\newtheorem{Cor}[Thm]{Corollary}
\newtheorem{Lem}[Thm]{Lemma}
\newtheorem{Prop}[Thm]{Proposition}
\newtheorem{Fac}[Thm]{Fact}

\theoremstyle{definition}
\newtheorem{Def}[Thm]{Definition}
\newtheorem{Rem}[Thm]{Remarks}

\theoremstyle{remark}

\def \cal{\mathcal}

\def\Ndb{ { {\rm I}\!{\rm N}} }
\def\Qdb{ {\rm Q}\kern-.65em {}^{ {}_/ }}
\def\Rdb{ {\rm I}\!{\rm R} }

\begin{document}
\centerline{SUBSPACES OF $c_0(\Ndb)$ AND LIPSCHITZ ISOMORPHISMS}
\medskip
\centerline{Gilles GODEFROY, Nigel KALTON, Gilles LANCIEN}
\vskip 1cm
\centerline{Abstract}

We show that the class of subspaces of $c_0(\Ndb)$ is
stable under Lipschitz isomorphisms. The main corollary is that
any Banach space which is Lipschitz-isomorphic to $c_0(\Ndb)$ is
linearly isomorphic to $c_0(\Ndb)$. The proof relies in part on
an isomorphic characterization of subspaces of
$c_0(\Ndb)$ as separable spaces having an equivalent norm such that the
weak-star and norm topologies quantitatively agree on the
dual unit sphere . Estimates on the Banach-Mazur distances are provided
when the Lipschitz constants of the isomorphisms are small.
The quite different non separable theory is also investigated.
\section{Introduction}

Banach spaces are usually considered within the category of
topological vector spaces, and isomorphisms between them are
assumed to be continuous and linear. It is however natural to
study them from different points of view, e.g. as infinite
dimensional smooth manifolds, metric spaces or uniform
spaces, and to investigate whether this actually leads to
different isomorphism classes. We refer to [J-L-S] and
references therein for recent results and description of this
field.  Some simply stated questions turn out to be hard to
answer: for instance, no examples are known of separable
Banach spaces $X$ and $Y$ which are Lipschitz isomorphic but not
linearly isomorphic. It is not even known if this could occur
when $X$ is isomorphic to $l_1$. The main result of this work is
that any separable space which is Lipschitz isomorphic to
$c_0(\Ndb)$ is linearly isomorphic to $c_0(\Ndb)$. Showing it
will require the use of various tools from non linear
functional analysis, such as the Gorelik principle. New
linear results on subspaces of $c_0(\Ndb)$ will
also be needed.

\medskip
 We now turn to a detailed description of our results. Section 2 contains
the main theorems of our article (Theorems 2.1 and 2.2),
which contribute to the classification of separable Banach spaces under
Lipschitz isomorphisms. These results are non linear.
However, their proof requires linear tools such as Theorem 2.4 which
provides a characterization of linear subspaces of
$c_0(\Ndb)$  in
terms of existence of equivalent norms with a property of
asymptotic uniform smoothness. This technical property is easier to handle
through the
dual norm, which is such that the weak* and norm topologies
agree quantitatively on the sphere (see Definition 2.3).
The main topological argument we need is
Gorelik's principle (Proposition 2.7) which is combined with a
renorming technique and with Theorem 2.4 for showing (Theorem
2.1) that the class of subspaces of
$c_0(\Ndb)$ is stable under Lipschitz-isomorphisms. It follows
(Theorem 2.2) that a Banach space is isomorphic to
$c_0(\Ndb)$ as soon as it is Lipschitz-isomorphic to it. The renorming
technique is somewhat similar to ``maximal rate of change"
arguments which are used for differentiating Lipschitz functions (see [P]).

We subsequently investigate extensions of the separable isomorphic results
of section 2 in two directions: what can be said when
the Lipschitz constants of the Lipschitz isomorphisms are small? What
happens in the non separable case? These questions are
answered in the last three sections. For reaching the answers, we have to
use specific tools, since the proofs are not
straightforward extensions of those from section 2.

Section 3 deals with  quantitative versions of Theorem 2.2. These
statements are ``nearly isometric" analogues, in the case of
$c_0(\Ndb)$, of Mazur's theorem which states that two isometric Banach
spaces are linearly isometric.
Indeed we show that a Banach space $X$ is close to
$c_0(\Ndb)$ in Banach-Mazur distance if there is a Lipschitz-isomorphism
$U$ between $X$ and
$c_0(\Ndb)$ such that the Lipschitz constants of $U$ and
$U^{-1}$ are close to 1 (Propositions 3.2 and 3.4). Proposition 3.2 relies
on an examination of the proof of Gorelik's principle in
the case of $c_0(\Ndb)$ and on an unpublished result of M. Zippin ([Z3]),
while Proposition 3.4 uses the concept of $K_0$-space
from [K-R].

The non separable theory is studied in sections 4 and 5. It is shown in
([J-L-S], Theorem 6.1) that if $1<p<\infty$, any Banach
space which is uniformly homeomorphic (in particular, Lipschitz isomorphic)
to $l^p(\Gamma)$ is linearly isomorphic to it, for any
set
$\Gamma$. But in the case of $c_o(\Gamma)$ (i.e. in the case $p=\infty$),
this situation happens to be quite different. Indeed there
are spaces which are Lipschitz isomorphic to
$c_0(\Gamma)$ with
$\Gamma$ uncountable but not linearly isomorphic to a subspace of that
space (see [D-G-Z 2] and Examples 4.9). The main gist of the
last two sections is that the separable theory extends to the class of
weakly compactly generated spaces (that is, to spaces $X$
which contain a weakly compact subset which spans a dense linear subspace)
but not further. Section 4 is devoted to characterizing
subspaces of $c_0(\Gamma)$ by the existence of certain equivalent norms,
that is, to extend Theorem 2.4 to the non-separable
case. It so happens that the quantitative behaviour of the equivalent
asymptotically uniformly smooth norms on
$X$, which  does not really matter in the separable case, is crucially
important in the non separable situation (Lemma 4.2).
Knowing this, we characterize, both isomorphically (Theorem 4.4) and almost
isometrically (Proposition 4.5) subspaces of
$c_0(\Gamma)$. We also obtain satisfactory classification results
for
$C(K)$-spaces, when some finite derivative of the compact space $K$ is
empty. More precisely, we show that
$K$ is an Eberlein compact and
$K^{(\omega_0)}$ is empty if and only if $C(K)$ is linearly
isomorphic to some space
$c_0(\Gamma)$ (Theorem 4.7), while the same equivalence holds with ``Lipschitz
isomorphic" if we drop the requirement that $K$ is Eberlein (Theorem 4.8).
Projectional resolutions of identity play a leading
role in this fourth section. The last section 5 contains the extension of
our main results Theorems 2.1 and 2.2 to the non
separable frame, which holds under the assumption of weakly compact
generation (Corollary 5.2). We also provide characterizations of
the spaces
$c_0(\Gamma)$, as well as some additional remarks about the non separable
results.
Our statements proved under the assumption ${\rm dens}(X)=\omega_1$
(Theorem 4.4, results from section 5) can be
extended with similar proofs to
the case ${\rm dens}(X) < \aleph_{\omega_0}$. It is plausible that cardinality
restrictions are not necessary.

A similar theory can be developed for uniform homeomorphisms.
This is the subject of the forthcoming paper [G-K-L2]. Some
results of [G-K-L2] and of the present paper have been
announced in [G-K-L1].

\medskip
{\sl Notation}: We denote by $B_X$, respectively $S_X$, the open unit
ball, respectively the
unit sphere of a Banach space $X$. If $V$ is a uniformly
continuous map from
a Banach space $X$ to a Banach space $Y$, we denote, for $t>0$,
$\omega(V,t)=
\sup \{\Vert Vx_1-Vx_2\Vert,\ \Vert x_1-x_2 \Vert \leq t \}$ its modulus
of
uniform continuity. Two Banach spaces $X$ and $Y$ are Lipschitz
isomorphic if there is a bijective map $U$ from $X$ onto $Y$
such that $U$ and $U^{-1}$ are both Lipschitz maps when $X$ and
$Y$ are equipped with the metric given by their norm. The word
``isomorphic", when used alone, will always mean linearly isomorphic.
The Lipschitz weak-star Kadec-Klee property (in short, $LKK^*$) is defined
in Definition 2.3, and in Definition 4.1 in the
non-separable case. We refer to the discussion that follows Definition 2.3
for the relation between this notion and V. Milman's
moduli from [M], and for related terminology. Specific notions which are
used in the non separable sections 4 and 5 are recalled
after Definition 4.1.

\medskip
{\sl Acknowledgement}: This work was initiated while the first and
last named authors were visiting the University of
Missouri-Columbia in 1997. They express their warmest thanks
to the Department of Mathematics of U. M. C. for its
hospitality and support, to their home institutions (C. N. R.
S., Universit\'e de Besan\c con) and to all those who made
this stay possible. The second author was supported in part by
NSF Grant DMS-9870027.

\section{SEPARABLE ISOMORPHIC RESULTS}

The aim of this section is to prove our main results.

\begin{Thm} The class of all Banach spaces that are linearly
isomorphic to a subspace of $c_0(\Ndb)$ is stable under Lipschitz
isomorphisms.\end{Thm}

When dealing  with $c_0(\Ndb)$ itself, we obtain a more precise theorem.

\begin{Thm} A Banach space is linearly isomorphic to $c_0(\Ndb)$ if
and only if it is  Lipschitz isomorphic to $c_0(\Ndb)$.\end{Thm}

The proof of these theorems will require the use of a characterization of
subspaces of $c_0(\Ndb)$ in terms of equivalent norms,
and a non-linear argument which relies mainly upon the Gorelik principle.
We first establish the somewhat technical renorming
characterization, and for this purpose we need to introduce some notation.
The following definition is consistent with the
terminology of ([K-O-S]); see also [D-G-K] and references therein.

\begin{Def} Let $X$ be a separable Banach space. The norm of
$X$  is said to be {\sl Lipschitz weak-star Kadec-Klee} (in short, $LKK^*$)
if there exists $c$ in
$(0,1]$ such
that  its dual norm satisfies the
following property: for any $x^*$ in $X^*$ and any weak$^*$ null
sequence
$(x^*_n)_{n \geq 1}$ in $X^*$
($x_n^* \buildrel {w^*}\over {\longrightarrow} 0$),
$$ \limsup \Vert x^*+x^*_n \Vert \geq \Vert x^*\Vert +
c\,\limsup \Vert x^*_n \Vert.$$
\noindent If the above property is satisfied with a given $c$ in
$(0,1]$, we
will say that the norm of $X$ is {\sl $c$-LKK$^*$}. If it is
satisfied with the optimal value
$c=1$,  we will say that the norm of  $X$ is {\sl
metric-KK$^*$}.\end{Def}

Let us make it clear that the above notion is a property of the norm of $X$
which is actually checked on the dual norm, the
reference to $X$ being contained in the use of the weak* topology. In this
paper, it will in practice be easier to work with dual
norms. However, it is appropriate to reformulate the above definition in
terms of the norm of $X$. This can be done using a modulus
which has been introduced in 1971 by V. Milman ([M]), and which we recall
now. If $x\in S_X$ and $Y$ is a linear subspace of $X$,
we let $$\overline{\rho}(\tau,x,Y)= sup\{\Vert x+y\Vert -1;~y\in Y,~\Vert
y\Vert\leq\tau\}$$ and then
$$\overline{\rho}(\tau,x)=inf\{\overline{\rho}(\tau,x,Y);~dim(X/Y)<\infty\}$$ an
d finally
$$\overline{\rho}(\tau)=sup\{\overline{\rho}(\tau,x);~x\in S_X\}$$ In
Milman's [M]
notation, $\overline{\rho}(\tau,x)=\delta(\tau;x,{\bf B}^{o})$.  In
[J-L-P-S], Banach spaces which satisfy that
$\overline{\rho}(\tau)=o(\tau)$ when $\tau$ tends to $0$ are called {\sl
asymptotically uniformly smooth}. An easy duality argument
shows that the norm of $X$ is Lipschitz weak-star Kadec-Klee if and only if
there exists $\tau_0>0$ such that
$\overline{\rho}(\tau_0)=0$. In fact, Lemma 2.5 below shows that if the
norm of $X$ is $c-LKK^*$ then $\overline{\rho}(c)=0$
and it follows from [K-W] that the norm of
$X$ is
$metric-KK^*$ if and only if
$\overline{\rho}(1)=0$. Hence, following the terminology of [J-L-P-S],
spaces which enjoy the $LKK^*$ property should be called
{\sl asymptotically uniformly flat}.  Although this latter terminology is
certainly more descriptive, we will keep using in the
statements the Lipschitz weak-star Kadec-Klee terminology since we
crucially use the dual presentation and the parameter $c$.

The following theorem asserts that having an equivalent
$LKK^*$ (or if preferred, asymptotically uniformly flat)
norm is an isomorphic characterization of the subspaces of $c_0(\Ndb)$.
The
precise  quantitative version of this result is the following.

\begin{Thm} Let $c$ in $(0,1]$ and $X$ be a separable Banach space
whose norm is
$c$-Lipschitz weak-star Kadec-Klee. Then, for any $\varepsilon >0$, there is a
subspace $E$
of $c_0(\Ndb)$ such that $d_{BM}(X,E) \leq {\displaystyle {1/
c^2}}+\varepsilon$; where
$d_{BM}(X,E)$ denotes the Banach-Mazur distance between $X$ and
$E$.\end{Thm}

\begin{proof} The following lemma gives a dual formulation of the notion
of
$LKK^*$ norm.

\begin{Lem} Let $c$ in $(0,1]$ and $X$ be a separable Banach space
with a $c$-Lipschitz weak-star Kadec-Klee norm, then
$$
{\rm max}\ (\Vert x \Vert,
{1 \over 2-c}\limsup \Vert x_n \Vert) \leq \limsup \Vert x+x_n \Vert
\leq  {\rm max}\ (\Vert x \Vert,{\displaystyle {1 \over c}}\,
\limsup \Vert x_n \Vert),
$$
whenever $(x_n)$ is a weakly null sequence in $X$
($x_n \buildrel {w}\over {\longrightarrow} 0$). Let us call
$m_\infty(c)$
this property.\end{Lem}

\begin{proof} Let $x$ in $X$ and $(x_n) \subset X$ with
$x_n \buildrel {w}\over {\longrightarrow} 0$. Without loss of
generality, we
may assume that $\lim\Vert x_n \Vert$ and $\lim\Vert x+x_n \Vert$ exist
with
$\lim\Vert x_n \Vert > 0$.

\noindent
We will first prove the right hand side inequality. For $n \geq 1$, pick
$y_n^*$
in $X^*$ so that $\Vert y_n^* \Vert =1$ and $y_n^*(x+x_n)= \Vert x+x_n
\Vert$.
Passing to a subsequence, we may assume that
$y^*_n \buildrel {w^*}\over {\longrightarrow} y^*$ and
$\lim\Vert y_n^*-y^* \Vert$ exists. Then, it follows from our assumption
that
$$
c \,\lim\Vert y_n^*-y^* \Vert \leq 1- \Vert y^* \Vert.
$$
Notice now that
$\Vert x+x_n \Vert=y^*_n(x+x_n)=y_n^*(x)+y^*(x_n)+(y_n^*-y^*)(x_n)$. So
we have
$$
\lim\Vert x+x_n \Vert \leq \Vert y^* \Vert \, \Vert x \Vert
+{(1-\Vert y^* \Vert)\, \lim\Vert x_n \Vert \over c} \leq
{\rm Max}\ (\Vert x \Vert,{\displaystyle {1 \over c}}\,
\lim\Vert x_n \Vert).
$$

\noindent
For the left hand side inequality, we only need to show that
$\lim\Vert x+x_n\Vert \geq {\displaystyle {1 \over 2-c}}\,
\lim\Vert x_n \Vert$. So we select now $x_n^*$ in $X^*$ with $\Vert
x_n^*
\Vert =1$ and $x_n^*(x_n)=1$ and we assume that
$x^*_n \buildrel {w^*}\over {\longrightarrow} x^*$ and
$\lim\Vert x^*_n -x^* \Vert$ exists. Again , we have
$$
c \,\lim\Vert x_n^*-x^* \Vert \leq 1- \Vert x^* \Vert.
$$
Since $(x_n^*-x^*)(x_n) \to \lim\Vert x_n \Vert$, we also obtain
$\lim\Vert x_n^*-x^* \Vert \geq 1$ and therefore $\Vert x^* \Vert \leq
1-c$.
We can write $x_n^*(x+x_n)=\Vert x_n \Vert + (x_n^*-x^*)(x) +x^*(x)$.
So, passing
to the limit we obtain $\lim\Vert x+x_n \Vert + (1-c) \Vert x \Vert \geq
\lim\Vert x_n \Vert$. Then we conclude by using the fact that
$\Vert x\Vert\leq \lim\Vert x+x_n \Vert$. \end{proof}

\noindent {\bf Remark.} The best constant ${\displaystyle {1/
(2-c)}}$ is
not crucial  for the proof of Theorem 2.2 that will be achieved with
the
trivial value
${\displaystyle {1/2}}$. However it will be used in the proof of
Proposition 3.2 and it helps us to relate this with  Theorem 3.2 in
[K-W]
which states, in the particular case
$p=\infty$, that a space  satisfying the property $m_\infty=m_\infty(1)$
embeds almost isometrically into
$c_0(\Ndb)$.

\medskip
Our next Lemma is the analogue of Lemma 3.1 in [K-W].

\begin{Lem}

\noindent
(i) If $F$ is a finite dimensional subspace of $X$ and $\eta>0$, then
there is a  finite dimensional subspace $U$ of $X^*$ such that
$$
\forall (x,y)\in F \times U_\perp \ \ (1-\eta)\,{\rm Max}
(\Vert x \Vert,{1 \over 2}\Vert y \Vert ) \leq \Vert x+y \Vert \leq
(1+\eta)\,{\rm Max}(\Vert x \Vert,{\displaystyle {1 \over c}}\Vert y
\Vert ).
$$
(ii) If $G$ is a finite dimensional subspace of $X^*$ and $\eta >0$,
then there is a
finite dimensional subspace $V$ of $X$ such that
$$
\forall (x^*,y^*)\in G \times V^\perp \ \ (1-\eta)
(\Vert x^* \Vert+c\Vert y^* \Vert ) \leq \Vert x^*+y^* \Vert \leq
\Vert x^* \Vert+\Vert y^* \Vert .
$$\end{Lem}

\begin{proof} Since the norm of $X$ is $LKK^*$, $X^*$ is
separable.
Then the proof is identical with the proof of Lemmma 3.1 in [K-W].
\end{proof}

We will now proceed with the proof of Theorem 2.4, which is only a slight
modification of the proof of Theorem 3.2 in [K-W]. So
let
$0<
\delta < {\displaystyle {1
\over 3}}$ and  pick  a positive integer $t$ such that $t>{\displaystyle
6(1+\delta) \over {c^3\delta}}$. Let also
$(\eta_n)_{n\geq 1}$  be a sequence of positive real numbers
satisfying
$$
0<\eta_n<{\displaystyle {\delta \over 2}},\ \
\prod_{n \geq 1}(1-\eta_n) >1-\delta\ {\rm and}\
\prod_{n \geq 1}(1+\eta_n) <1+\delta.
$$
Finally, let $(u_n)_{n\geq 1}$ be a dense sequence in $X$. Following the
ideas of Kalton and
Werner, we then construct subspaces $(F_n)_{n \geq 1}$, $(F'_n)_{n \geq
1}$ of
$X^*$ and
 $(E(m,n))_{1 \leq m \leq n}$ of $X$ so that:

\smallskip\noindent
(a) dim $F_n < \infty$, dim $E(m,n)<\infty$ for all $m \leq n$.

\smallskip\noindent
(b) $F'_n\subseteq [u_1,..,u_n]^\perp \cap \bigcap_{j \leq
k<n}E(j,k)^\perp$
is weak$^*$-closed and $X^*=F_1\oplus..\oplus F_n  \oplus F'_n$.

\smallskip\noindent
 (c) $F'_n=F_{n+1} \oplus F'_{n+1}$.

\smallskip\noindent
(d) If $x^* \in F_1+..+F_n$ and $y^* \in F'_{n+1}$, then
$$(1-\eta_n)(\Vert x^* \Vert+c\Vert y^* \Vert ) \leq \Vert x^*+y^* \Vert
\leq
\Vert x^* \Vert+\Vert y^* \Vert.$$

\smallskip\noindent
(e) If $x \in (F_1+..+F_n)_\perp$ and $y \in \sum_{j \leq k <n}E(j,k)$,
then
$$
(1-\eta_n)\,{\rm Max}(\Vert x \Vert,{1 \over 2}\Vert y \Vert ) \leq
\Vert x+y \Vert \leq
(1+\eta_n)\,{\rm Max}(\Vert x \Vert,{\displaystyle {1 \over c}}\Vert y
\Vert ).
$$

\smallskip\noindent
(f) $(F_1+..+F_{m-1}+F'_n)_\perp \subseteq E(m,n)$ and $E(m,n)\subseteq
(F_1+..+F_{m-2})_{\perp}$ if $1 \leq m \leq n$.

\smallskip\noindent
(g) If $x^* \in F_m+..+F_n$, then there exists $x \in E(m,n)$ so that
$\Vert x \Vert \leq 1$ and

$x^*(x) \geq c(1-\delta)\Vert x^*\Vert$.

\smallskip\noindent
Now, as in [K-W], we define, for $0\leq s\leq t-1$
$$
T_s:Y_s=c_0(E(4(n-1)t+4s+4,4nt+4s+1)_{n \geq 0}) \to X
$$
and
$$
R_s:Z_s=c_0(E(4nt+4s+2,4nt+4s+3)_{n \geq 0}) \to X
$$
by $T_s((y_n)_{n \geq 0})=\sum y_n$ and $R_s((z_n)_{n \geq 0})=\sum
z_n$. And
also
$$
T:Y=\ell_\infty((Y_s)_{s=0}^{t-1}) \to X\ \ {\rm and}\ \
R:Z=\ell_\infty((Z_s)_{s=0}^{t-1}) \to X
$$
by
$$
T(\xi_0,..,\xi_{t-1})={\displaystyle {1 \over t}}\sum_{s=0}^{t-1}
T_s\xi_s\ \
{\rm and}\ \ R(\xi_0,..,\xi_{t-1})=\sum_{s=0}^{t-1} R_s\xi_s.
$$
Then we get
$$
\forall \xi \in Y_s,\ {1-\delta \over 2}\Vert \xi \Vert \leq \Vert T_s
\xi
\Vert \leq  {\displaystyle {1+\delta \over c}}\Vert \xi \Vert,
$$
$$
\forall \xi \in Z_s,\ {1-\delta \over 2}\Vert \xi \Vert \leq \Vert R_s
\xi
\Vert \leq  {\displaystyle {1+\delta \over c}}\Vert \xi \Vert,
$$
$$
\Vert T \Vert \leq {1+\delta \over c}\ \ {\rm and}\ \
\Vert R \Vert \leq {1+\delta \over c}.
$$
Still following [K-W] we can also show that if $x^*$ in $X^*$ satisfies
$R_s^*x^*=0$, then

\noindent
$\Vert T_s^*x^* \Vert \geq c(1-\delta)\Vert x^* \Vert$. Then a
Hahn-Banach argument yields
$$
\forall x^* \in X^*\ \Vert T_s^* x^*\Vert \geq c(1-\delta)\Vert x^*\Vert
-
2(c+{\displaystyle {1+\delta \over c(1-\delta)}})\Vert R_s^*x^* \Vert.
$$
Since $\delta < {\displaystyle {1 \over 3}}$ and $c \leq 1$, we have
$$\forall x^* \in X^*\ \Vert T_s^* x^*\Vert \geq
c(1-\delta)\Vert x^*\Vert -{\displaystyle {6 \over c}}\Vert
R_s^*x^* \Vert.$$ Therefore
$$\forall x^* \in X^*\ \Vert T^* x^*\Vert \geq c(1-\delta)\Vert
x^*\Vert -
{6 \over ct}\Vert R^*x^* \Vert \geq \Big\lbrack c(1-\delta) -
{6(1+\delta) \over c^2t}\Big\rbrack \,
\Vert x^*\Vert.$$
Thus, for our initial choice of $t$ we obtain
$$\forall x^* \in X^*\ \Vert T^*x^* \Vert \geq c(1-
2\delta)\Vert x^*\Vert.$$ Since we have on the other hand that
$\Vert T \Vert \leq {\displaystyle {1 \over c}}(1+\delta)$, we
get that
$$
d(X,Y/{\rm ker}T) < {1+\delta \over c^2(1-2\delta)}.
$$
As a
$c_0$-sum of finite dimensional spaces, $Y$ embeds almost isometrically
into
$c_0(\Ndb)$. Then, by Alspach's theorem [Al], so does $Y/{\rm ker}T$.
This concludes our proof. \end{proof}

Let us mention that Theorem 2.4 is much easier to show through a skipped
blocking argument when the space $X$ is assumed to have a
shrinking
$FDD$. Now, using [J-R] and the simple fact that being a subspace of
$c_0(\Ndb)$ is a three-space property, the general case
follows. This alternative approach from [J-L-P-S] does not provide however
the same isomorphism constants.

We now turn to non linear theory. First we state a slight modification of
the Gorelik Principle as
it is  presented in [J-L-S].

\begin{Prop} {\rm (Gorelik's Principle)} -  Let $E$ and $X$ be two
Banach spaces  and $U$ be a homeomorphism from $E$ onto $X$ with
uniformly
continuous inverse.  Let $b$ and $d$ two positive constants and let
$E_0$ be a
subspace of finite  codimension of $E$. If $d> \omega(U^{-1},b)$, then
there
exists a compact subset
$K$ of $X$ such that
$$ bB_X \subset
K+U(2dB_{E_0}).$$\end{Prop}

\begin{proof} We recall a fundamental lemma due to E.
Gorelik
[G] and that can also be found in [J-L-S].

\begin{Lem} For every $\varepsilon >0$ and $d>0$, there exists a
compact subset $A$ of
$dB_E$ such that, whenever $\Phi$ is a continuous map  from
$A$ to $E$ satisfying $\Vert \Phi (a) -a \Vert < (1-\varepsilon)d$ for
any
$a$ in $A$, then $\Phi(A) \cap E_0 \neq \emptyset$.\end{Lem}

 Now, fix $\varepsilon >0$ such that
$d(1-\varepsilon)>\omega(U^{-1},b)$. Let $K=-U(A)$, where
$A$ is the compact set obtained in Lemma 2.8. Consider now
$x$ in $ bB_X$ and the map
 $\Phi$ from $A$ to $E$ defined by $\Phi(a)=U^{-1}(x+Ua)$. It is clear
that
for any $a$ in $A$,
$\Vert \Phi (a) -a \Vert < (1-\varepsilon) d$. Then, it follows from
Lemma 2.8
that there exists
$a \in A$ so that $U^{-1}(x+Ua) \in 2dB_{E_0}$. This concludes our
proof.
\end{proof}

We can now proceed to prove Theorem 2.1.

\begin{proof}{\sl of Theorem 2.1}:   Let $U$ be a Lipschitz isomorphism
from a subspace $E$ of
$c_0$
onto the Banach space $X$. Theorem 2.4 indicates that we need to build
an
equivalent $LKK^*$ norm on $X$. This norm will be defined as
follows.
For $x^*$ in $X^*$, set:
$$ |||x^*|||= \sup \{ {|x^*(Ue-Ue')| \over ||e-e'||};\ (e,e') \in E
\times E,\
e \neq e'\}.$$
Since $U$ and  $U^{-1}$ are Lipschitz maps, $|||\ |||$ is an equivalent
norm
on $X^*$. It is clearly weak$^*$ lower semicontinuous and therefore is
the
dual norm of an equivalent norm on $X$ that we will also denote
$|||\ |||$.

Consider $\varepsilon >0$, $x^* \in X^*$ and $(x^*_k)_{k\geq 1} \subset
X^*$
such that $x_k^* \buildrel {w^*}\over {\longrightarrow} 0$ and

\noindent $||x^*_k||\geq\varepsilon>0$ for all $k \geq 1$.  Fix $\delta
>0$ and
then
$e$ and $e'$ in $E$ so that
$${\displaystyle{x^*(Ue-Ue') \over ||e-e'||}} >(1-\delta)|||x^*|||.$$
By using
translations in order to modify $U$, we may as well assume that
$e=-e'$ and $Ue=-Ue'$. Since $E$ is a subspace of $c_0$, it admits a
finite
codimensional subspace $E_0$ such that
$$\forall f \in
||e|| B_{E_0},\ \  ||e+f|| \vee ||e-f|| \leq
(1+\delta)||e||.\eqno(2.1)$$
Let $C$ be the Lipschitz constant of $U^{-1}$. By Proposition 2.7, for
every
$b<\displaystyle{||e|| \over 2C}$ there is a compact subset
$K$ of $X$ such that $bB_X \subset K +U(||e||B_{E_0})$. Since $(x_k^*)$
converges uniformly to $0$ on any compact subset of $X$, we can
construct a
sequence $(f_k) \subset ||e||B_{E_0}$ such that:
$$\liminf x^*_k(-Uf_k) \geq {\varepsilon ||e|| \over 2C}.$$

\noindent We deduce from (2.1) that $x^*(Uf_k+Ue) \leq
(1+\delta)||e||\,|||x^*|||$ and therefore $x^*(Uf_k) \leq 2\delta
||e||\,|||x^*|||$. Using again the fact that
$x_k^* \buildrel {w^*}\over {\longrightarrow} 0$, we get that:
$$\liminf (x^*+x^*_k)(Ue-Uf_k) \geq (1-3\delta)||e||\,|||x^*||| +
{\varepsilon ||e|| \over 2C}.$$
Since $\delta$ is arbitrary, by using the definition of $|||\ |||$
and (2.1),  we obtain
$\liminf |||x^*+x^*_k||| \geq |||x^*||| + {\displaystyle{\varepsilon
\over
4C}}$. This proves that $|||\ |||$ is $LKK^*$, and concludes the proof of
Theorem 2.1.\end{proof}

Theorem 2.2 is easily deduced from Theorem 2.1 through the use of two
classical results.

\begin{proof}{\sl of Theorem 2.2}: We only need to prove the ``if" part. So
let $X$ be a
Banach
space which is Lipschitz isomorphic to $c_0(\Ndb)$. Theorem 2.1 asserts
that
$X$ is linearly isomorphic to a subspace of $c_0(\Ndb)$. Besides, it is
known
that the class of all ${\cal L}^\infty$ spaces is stable under uniform
homeomorphisms ([H-M]) and that a ${\cal L}^\infty$ subspace of
$c_0(\Ndb)$ is
isomorphic to $c_0(\Ndb)$ ([J-Z]). This establishes Theorem 2.2.\end{proof}

Note that although Theorems 2.1 and 2.2 are non linear results, it is
critically important that the Banach space $X$ is Lipschitz
isomorphic to a {\sl linear subspace} of $c_0(\Ndb)$. In fact, given any
separable Banach space $Y$, there is a bi-Lipschitz map
between
$Y$ and a {\sl subset} of $c_0(\Ndb)$ ([Ah]).

\section{QUANTITATIVE RESULTS}

Recall that for $\lambda \geq 1$, a Banach space $X$ is said to be
${\cal
L}^\infty_\lambda$ if for every finite dimensional subspace $E$ of $X$,
there
is a finite dimensional subspace $F$ of $X$, containing $E$ and such
that
$d_{BM}(F,\ell_\infty^{{\rm dim}\,F})\leq \lambda$. If $X$ is ${\cal
L}^\infty_\lambda$ for some $\lambda \in [1,+\infty)$, then it is said
to be
${\cal L}^\infty$ (see [L-T]). We already used the fact ([J-Z]) that a
subspace of $c_0(\Ndb)$ is isomorphic to $c_0(\Ndb)$ if
and only
if it is
${\cal L}^\infty$. Combining this with
Theorem 2.4, we get that if a separable ${\cal L}^\infty$ space admits a
$LKK^*$ norm, then it is isomorphic to $c_0(\Ndb)$.
The following
statement gives a quantitative estimate on the linear isomorphism.

\begin{Prop} There exists a function $F:[1,+\infty)\times (0,1] \to
[1,+\infty)$ such that if
$X$ is a separable ${\cal L}^\infty_\lambda$ space with a
$c$-Lipschitz weak-star Kadec-Klee norm, then
$$d_{BM}(X,c_0(\Ndb))
\leq F(\lambda,c).$$ Moreover $F(1,1)=1$ and $F$ is continuous
at $(1,1)$.
\end{Prop}

\begin{proof} Let us first mention that for  values of $\lambda$ and $c$
close to $1$, the result follows directly from a work of M.
Zippin [Z3], who
proved that if $X$ is a ${\cal L}^\infty_\mu$ subspace of $c_0$ with
$\mu <7/6$, then
$$d_{BM}(X,c_0(\Ndb)) \leq {\mu^2 \over \mu^2 -2\mu^3+2}.$$
Then, it is easily checked that in our setting, if we assume moreover
that
$\lambda/c^2<7/6$, we get
$$d_{BM}(X,c_0(\Ndb)) \leq {\lambda^2 \over
2c^6+\lambda^2c^2-2\lambda^3}.$$

For the general case we do not have an explicit function $F$. We will
just
reproduce an argument by contradiction used in ([G-L], p.261).
Indeed, if there is no such function, then there exist
$(\lambda,c)$ in
$[1,+\infty)\times (0,1]$ and a sequence $(X_n)$ of separable ${\cal
L}^\infty_\lambda$ spaces with a
$c-LKK^*$ norm such that, for all $n\geq 1$,
$d_{BM}(X_n,c_0(\Ndb)) \geq n$. But the space
$Y=\big(\sum\oplus X_n\big)_{c_0}$ is ${\cal L}^\infty$ with a
$LKK^*$ norm and thus by Theorem 2.4 and [J-Z] it is isomorphic to
$c_0(\Ndb)$. So,
the
$X_n$'s being uniformly complemented in $Y$, their Banach-Mazur distance
to
$c_0(\Ndb)$ should be bounded, a contradiction.
\end{proof}

We will now give two quantitative versions of Theorem 2.2.

\begin{Prop} There exists a function $F:(1,+\infty) \to (1,+\infty)$
such that
$\displaystyle{\lim_{\lambda \to 1^+}} F(\lambda)= 1$ and such that if
$U:X \to
c_0(\Ndb)$ is a bi-Lipschitz map with ${\rm Lip}(U)\cdot{\rm
Lip(U^{-1})}=\lambda$,
then
$d_{BM}(X,c_0(\Ndb)) \leq F(\lambda)$.\end{Prop}

\begin{proof} Let $E_n=\{x=(x(i))_{i\geq 0} \in c_0(\Ndb);\ x(i)=0 \
{\rm if}\
i>n
\}$. We set $A=B_{E_n}$. It is easily seen that if $\Phi:A \to
c_0(\Ndb)$ is a
continuous map such that $\Vert a-\Phi(a)\Vert \leq 1$ for all $a \in
A$, then
there exists $a_0 \in A$ such that $\Phi(a_0)(i)=0$ for all $i \leq n$.
Indeed,
if $\pi: c_0(\Ndb) \to E_n$ is the natural projection and
$F(a)=a-\pi(\Phi(a))$, then $F(A) \subseteq A$ and by Brouwer's theorem,
there
is $a_0 \in A$ with
$F(a_0)=a_0$. Hence $\vert\Phi(a_0)(j)\vert \leq 1$ for all $j>n$, and
thus
$\Phi(a_0) \in B_{F_n}$, where
$$F_n=\{x \in c_0;\ x(j)=0 \ {\rm if}\ j \leq n\}.$$
If we now reproduce the proof of Gorelik's Principle (Proposition 2.7),
using
the compact set $A$ and the space $F_n$ defined above (with an
appropriate
choice of $n$), we find in the notation of the proof of Theorem 2.1
that for
any
$b<\displaystyle{(\Vert e \Vert /{\rm Lip}(U^{-1}))}$, there is a
compact
subset $K$ of $X$ such that
$$bB_X \subset K+U(\Vert e \Vert B_{F_n})$$
and it follows that the norm $|||~.~ |||$ is
$\lambda^{-1}-LKK^*$.
Now Theorem 2.4 shows that the distance from $(X,|||~.~ |||)$
to the subspaces
of $c_0(\Ndb)$ is at most $\lambda^2$. Since the distance between the
original
norm
$\Vert~.~ \Vert$ of $X$ and $|||~.~ |||$ is less than $\lambda$,
it follows that
the Banach-Mazur distance from $(X,\Vert~.~ \Vert)$ to the
subspaces of
$c_0(\Ndb)$ is at most $\lambda^3$.

\medskip
We now observe the following

\begin{Fac} There is a function $F_0:(1,+\infty) \to (1,+\infty)$ with
$\displaystyle{\lim_{\lambda \to 1^+}} F_0(\lambda)= 1$, and such that
if $X$
satisfies the assumptions of the proposition, then $X$ is an ${\cal
L}^\infty_{F_0(\lambda)}$ space.\end{Fac}

\begin{proof} By the ultrapower version of the local reflexivity
principle,
$X^{**}$ is isometric to a 1-complemented subspace of some ultrapower
$(X)_{\cal U}$. We set $Z=(c_0)_{\cal U}$. Clearly, there is a
bi-Lipschitz
map $\tilde U: (X)_{\cal U} \to Z$ with

\noindent Lip$(\tilde U)~\cdot$ Lip$(\tilde
U^{-1})=\lambda$. It follows that there are maps $f:X^{**} \to Z$ and
$g:Z \to
X^{**}$ with Lip$(f)~\cdot$ Lip$(g)=\lambda$ and $g \circ f =Id_{X^{**}}$. By
([H-M], Lemma 2.11.), there is $\tilde g:Z^{**} \to X^{**}$ extending
$g$ and
such that Lip$(\tilde g)$=Lip$(g)$. The space $Z^{**}$ is isometric to
the
dual of an $L^1$-space, hence it is a ${\cal P}_1$ space (see [L-T], p.162).
Since
$\tilde g \circ f = Id_{X^{**}}$, it follows that if $M$ is a metric
space,
$N$ a subspace of $M$ and $\psi:N \to X^{**}$ a Lipschitz map, there
exists a
Lipschitz extension $\overline \psi :M \to X^{**}$ with Lip$(\overline
\psi)
\leq \lambda\,$Lip$(\psi)$. In particular, $X^{**}$ is isometric to a
linear
subspace $Y$ of $l_\infty(\Gamma)$ on which there exists a Lipschitz
projection $P$ with Lip$(P) \leq \lambda$. Since $X^{**}$ is
1-complemented
in its own bidual, it follows from ([Li], Corollary 2 to Theorem 3) that
there
exists a linear projection $\pi:l_\infty(\Gamma) \to Y$ with $\Vert \pi
\Vert
\leq \lambda$. Therefore, $X^{**}$ is a ${\cal P}_\lambda$ space.

By ([L-R], see p. 338), $X^{**}$ is therefore a ${\cal
L}^\infty_{10\lambda}$ space, and so is $X$. Moreover ([Z1],[Z2] and [B]
Th.
13), when $F$ is a finite dimensional ${\cal P}_{1+\varepsilon}$ space
with $\varepsilon < 17^{-8}$, then if we let $\nu=\varepsilon^{1/8}$
(see [B])
$$d_{BM}(F,l_\infty^{{\rm dim}(F)}) \leq {1+6\nu \over
(1-6\nu)(1-17\nu)}.$$
In the above notation, any finite dimensional subspace of $Y$ is
contained, up
to $\delta >0$ arbitrary, in a space $\pi(G)$, where $G$ is isometric
to a
finite dimensional $l_\infty$. Such an $F$ is ${\cal P}_\lambda$;
therefore
([B], Theorem 13) guarantees the requirement
$\displaystyle{\lim_{\lambda \to
1^+}}F(\lambda)= 1$. \end{proof}

We now proceed with the proof of Proposition 3.2. We know that
$(X,\Vert ~.~
\Vert)$
is a ${\cal
L}^\infty_{F_0(\lambda)}$ space whose Banach-Mazur distance to the
subspaces of
$c_0(\Ndb)$ is at most $\lambda^3$. Any ${\cal
L}^\infty_{\mu}$ subspace $H$ of $c_0(\Ndb)$ is isomorphic to
$c_0(\Ndb)$
([J-Z]) and by using contradiction (see [G-L], p.261) we show the
existence of
a function
$F_1(\mu)$ such that $d_{BM}(H,c_0)\leq F_1(\mu)$. Finally according to
([Z3]), there is such a function $F_1$ which satisfies
$\displaystyle{\lim_{\mu \to
1^+}}F_1(\mu)= 1$ (see proof of Proposition 3.1 above). For any
$\varepsilon >0$,
$X$ is
$(\lambda^3+\varepsilon)$ isomorphic to a subspace $G_\varepsilon$ of
$c_0(\Ndb)$ which is a ${\cal L}^\infty_{\mu}$ space with
$\mu=(\lambda^3+\varepsilon)F_0(\lambda)$; the existence of $F$ as
claimed in
the proposition clearly follows. \end{proof}

Using the techniques from [K-O] and the notion of $K_0$-space ([K-R]),
we can
actually extend Proposition 3.2. to arbitrary equivalent renormings of
$c_0(\Ndb)$.

\begin{Prop} Let $Y$ be a Banach space which is linearly isomorphic to
$c_0(\Ndb)$. Then for any $\varepsilon >0$, there is $\delta >0$ such
that if
$X$ is a Banach space and $U:X \to Y$ is a bi-Lipschitz onto map with
Lip$(U)~\cdot$ Lip$(U^{-1})<1+\delta$, then $d_{BM}(X,Y)<1+\varepsilon.$
\end{Prop}

\begin{proof} The proof relies heavily on [K-O], from which we take the
following notation: if $d_M(E,F)$ denotes the Hausdorff distance between
two
subsets $E$ and $F$ of a metric space $M$, the Kadets distance
$d_K(X,Y)$
between two Banach spaces $X$ and $Y$ is
$$d_K(X,Y)= \inf\{d_Z(U(B_X),V(B_Y))\},$$
where the infimum is taken over all linear isometric embeddings $U,V$ of
$X,Y$
into an arbitrary common Banach space $Z$.

The Gromov-Hausdorff distance $d_{GH}(X,Y)$ is the infimum of the
Hausdorff
distances $d_M(B_X,B_Y)$ over all isometric embeddings of $X$
and $Y$ into an
arbitrary common metric space $M$. By ([K-O], Th. 2.1), we have
$$d_{GH}(X,Y)\leq \sup\{{1 \over 2}\Vert \Phi(x)-\Phi(x')\Vert_Y -\Vert
x-x'\Vert_X;\ x,x' \in B_X\}$$
where $\Phi:B_X \to B_Y$ is a bijective map.

It follows easily that for any $\eta >0$, there is $\delta >0$ such that
if
there is $U:X \to Y$ a Lipschitz isomorphism with
Lip$(U)~\cdot$ Lip$(U^{-1})<1+\delta$, then $d_{GH}(X,Y) < \eta$. Obviously,
one has
$d_{GH}(X,Y) \leq d_K(X,Y)$ (and in general these two distances are not
equivalent: for instance ([K-O]), $\displaystyle{\lim_{p \to
1^+}}\,d_{GH}(\ell_p,\ell_1)=0$ while $d_K(\ell_p,\ell_1)=1$ for all
$p>1$). We
recall that a Banach space $E$ is a $K_0$-space ([K-R]) if there exists
$K_0>0$ such that whenever $f:E \to \Rdb$ is a homogeneous function
which is
bounded on
$B_E$ and satisfies
$$\vert f(x+x')-f(x)-f(x')\vert \leq \Vert x\Vert +\Vert x'\Vert$$
then there exists $x^* \in E^*$ such that
$$\vert f(x)-x^*(x) \vert \leq K_0 \Vert x \Vert,\ \forall x \in E.$$
It is shown in [K-R] that $c_0(\Ndb)$ is a $K_0$-space. By ([K-O],
Theorem
3.7), if
$E$ is a $K_0$-space and $(E_n)$ is such that $\lim d_{GH}(E_n,E)=0$,
then
$\lim d_K(E_n,E)=0$.

Since $Y$ is a $K_0$-space as isomorphic to $c_0(\Ndb)$, for any $\alpha
>0$
there is $\eta >0$ such that $d_{GH}(X,Y)<\eta$ implies
$d_K(X,Y)<\alpha$. Let
$Z$ be a Banach space which contains isometric copies of $X$ and $Y$
with $d_Z(B_X,B_Y)<\alpha$. We may and do assume that $Z$ is separable.
By
Sobczyk's theorem, $Y$ is linearly complemented in any separable
super-space
$Z$, and the norm of the projection $\pi_Z$ is bounded independently of
$Z$.
It easily follows that given $\varepsilon >0$, there is $\alpha>0$ such
that
if $d_Z(B_X,B_Y)<\alpha$ then $d_{BM}(X,Y)<1+\varepsilon$. Indeed, the
restriction to $X$ of $\pi_Z$ provides the required linear isomorphism.
This
concludes the proof. \end{proof}

\section{SUBSPACES OF $c_0(\Gamma)$}
\medskip

We now consider non separable spaces. It turns out that the non
separable
theory looks quite different. In this section we first establish non
separable analogues of Theorem 2.4 for characterizing subspaces
of $c_0(\Gamma)$ spaces, then we determine which compact spaces $K$ are
such that the Banach space $C(K)$ is linearly or Lipschitz
isomorphic to a
$c_0(\Gamma)$ space. It turns out that the two properties are distinct in
the non separable case, and this leads to a bunch
of natural non separable spaces which are Lipschitz but not linearly
isomorphic to  $c_0(\Gamma)$.

As will be clear in the sequel, the techniques that we develop are separably
determined. So we adopt the following definition:

\begin{Def} Let $X$ be a Banach space, and let $c\in (0,1]$. The norm
$\Vert \ \Vert$ of $X$
is
$c$-Lipschitz weak-star Kadec-Klee if its restriction to any separable subspace
of
$X$ is $c$-Lipschitz weak-star Kadec-Klee in the sense of Definition 2.3.
\medskip
\noindent If $c=1$, we say again that the norm is metric weak-star Kadec-Klee.
\end{Def}

We now recall classical terminology from non separable Banach space theory.
A projectional resolution of identity (in short, P.R.I.)
is a well-ordered sequence of norm-one projections which allows to ``break"
a non separable Banach space into smaller subspaces. We
refer to ([D-G-Z], Chapter VI) or ([Di]) for a precise definition and basic
properties of projectional resolutions of identity . A
projectional resolution of identity
$(P_\alpha)$ is said to be shrinking when $(P_\alpha^*)$ is a P.R.I on
$X^*$. A Banach space $X$ is weakly compactly generated (in short, w.c.g.)
if it contains a weakly compact subset which spans a
dense linear subspace. By [A-L], every w.c.g. space has a projectional
resolution of identity.

 We now state and prove two lemmas which lead to our non separable analogue
of Theorem 2.4.

\begin{Lem} Let $X$ be a Banach space. If $\Vert \ \Vert$ is a
metric weak-star Kadec-Klee
norm on $X$, then  $(X,\Vert \ \Vert)$ has a shrinking projectional
resolution of identity, and thus
$X$ is
weakly compactly generated.
\end{Lem}

\begin{proof} By ([F-G], Th. 3), it suffices to show that $(X,\Vert \
\Vert)$
is an $M$-ideal in its bidual. But ([H-W-W], Cor.III.1.10), asserts
that to
be an $M$-ideal in its bidual is a separably determined property. So let
$Y$ be a separable subspace of $X$ and
$\pi: Y^{***}
\to Y^*$ be the canonical projection. Pick $t \in Y^{***}$
with $\Vert t \Vert=1$, and write $t=y^*+s$ with
$s \in {\rm Ker}\, \pi= Y^\perp$. By definition of an $M$-ideal, what we
need to show is that
$$\Vert t\Vert =\Vert y^*\Vert+\Vert s\Vert \eqno (4.1)$$
There is a net $(y^*_\alpha)$ in $B_{Y^*}$ such that $t= \lim
y^*_\alpha$ in
$(Y^{***},w^*)$ and then

\noindent $y^*=\lim y^*_\alpha$ in
$(Y^{*},w^*)$. Since $\Vert t-y^* \Vert =\Vert s \Vert$, we have
$$\lim\inf \Vert y^*_\alpha -y^* \Vert \geq \Vert s \Vert$$
and since $\Vert\ \Vert$ is metric-$KK^*$, this implies that $\Vert
y^*\Vert
\leq 1-\Vert s \Vert =\Vert t \Vert -\Vert s \Vert$. This shows (4.1)
since
$\Vert t \Vert \leq \Vert y^* \Vert+\Vert s \Vert$ by the triangle
inequality.
\end{proof}

\begin{Lem} Let $X$ be a Banach space with a $c$-Lipschitz weak-star
Kadec-Klee norm.
For
every $x \in X$, there exists a separable subspace $E$ of $X^*$ such
that if

\noindent$y\in E_\perp \subset X$, one has
$$ \max (\Vert x \Vert, {\Vert y \Vert \over 2-c}) \leq \Vert x+y \Vert
\leq \max(\Vert x \Vert, {\Vert y\Vert \over c}).$$
\end{Lem}

\begin{proof} It clearly suffices to show that for any $\varepsilon >0$,
there is $F \subset X^*$ separable such that if $y \in F_\perp$
$$ (1-\varepsilon)\max (\Vert x \Vert, {\Vert y \Vert \over 2-c}) \leq
\Vert
x+y
\Vert
\leq (1+\varepsilon)\max(\Vert x \Vert, {\Vert y\Vert \over c}).
$$
Assume, for instance, that for any separable $F$, there is $y\in
F_\perp$ such that
$$ \Vert
x+y
\Vert
> (1+\varepsilon)\max(\Vert x \Vert, {\Vert y\Vert \over c}).
$$
We construct inductively an increasing sequence $(F_n)$ of separable
subspaces
of $X^*$, and $(y_n)$ in $X$ such that for all $n \geq 1$,

\noindent (i) if $u \in {\rm span}\,\{x,y_1,..,y_n\}$, then
$\Vert u\Vert =\sup\{\vert f(u)\vert;\ \Vert f \Vert \leq 1,\ f \in
F_n\}$.

\noindent (ii) $y_{n+1} \in (F_n)_\perp$.

\noindent (iii) $\Vert x+y_{n+1} \Vert > (1+\varepsilon) \max(\Vert
x\Vert,\displaystyle{\Vert y_{n+1}\Vert \over c})$.

\noindent We let $G=\overline{\rm span} \{x,(y_j)_{j\geq 1}\}$. Since
the
weak$^*$ and norm topologies coincide on $S_{G^*}$, it follows from (i)
that
$D=\cup_{n \geq 1} (F_n)_{|_G}$ is dense in $(G^*,\Vert\  \Vert)$. Then
(ii)
implies that $y_n \buildrel {w}\over {\longrightarrow} 0$. But now (iii)
contradicts Lemma 2.5. This proves the lemma, since we can
clearly proceed
along the same lines with the left hand side of the inequality.
\end{proof}

We now state and prove an analogue to Theorem 2.4 for non
separable spaces. To avoid dealing with singular cardinals, we limit
ourselves
to the case where the density character of $X$, denoted by dens$(X)$, is
equal
to $\omega_1$. It is plausible that this restriction is irrelevant.

\begin{Thm} Let $X$ be a Banach space such that {\rm
dens}$(X)=\omega_1$.
Then
$X$ is weakly compactly generated and $X$ has an equivalent Lipschitz
weak-star Kadec-Klee norm if and
only if
$X$ is isomorphic to a subspace of $c_0(\Gamma)$, where $\vert \Gamma
\vert=\omega_1$.\end{Thm}

\begin{proof} The natural norm of $c_0(\Gamma)$ is metric-$KK^*$ and
every
subspace of $c_0(\Gamma)$ is w.c.g. ([Jo-Z], see also [D-G-Z], Chapter
VI).

Conversely, if $X$ is w.c.g. and has a $c-LKK^*$ norm, then
$X$ is
a w.c.g. Asplund space and thus ([F], see also [D-G-Z], Th VI.4.3) $X$
has a
shrinking P.R.I. $(P_\alpha)_{\alpha \leq \omega_1}$. Using Lemma 4.3,
we
construct by induction on $\alpha$, ordinals $\lambda_\alpha<\omega_1$
such
that $\lambda_\alpha<\lambda_{\alpha+1}$ and such that if
$P_{\lambda_\alpha}(x)=x$ and $P_{\lambda_\alpha}(y)=0$, then
$$ \max(\Vert x\Vert,{\Vert y\Vert \over 2-c}) \leq \Vert x+y\Vert
\leq \max(\Vert x\Vert ,{\Vert y\Vert \over c}).
$$
If we let $X_\alpha=(P_{\lambda_{\alpha+1}}-P_{\lambda_\alpha})(X)$,
then $X$
is isomorphic to $(\sum \oplus X_\alpha)_{c_0}$. By Theorem 2.4, the
spaces
$X_\alpha$ are (uniformly in $\alpha$) isomorphic to subspaces of
$c_0(\Ndb)$;
this concludes the proof.
\end{proof}

We now provide a nearly isometric result. It follows from Lemma 4.2 and
Theorem 4.4 that any space $X$ with
dens$(X)=\omega_1$ which has a metric-$KK^*$ norm is isomorphic to a
subspace
of $c_0(\Gamma)$ with $|\Gamma|=\omega_1$. However, a much better result
is
available, namely:

\begin{Prop} Let $X$ be a Banach space. The following assertions are
equivalent:

\noindent (i) The norm of $X$ is metric weak-star Kadec-Klee.

\noindent (ii) For any $\varepsilon>0$, there
is a subspace
$X_\varepsilon$ of
$c_0(\Gamma)$, with
$|\Gamma|={\rm dens} (X)$, such that
$d_{BM}(X,X_\varepsilon)<1+\varepsilon$.\end{Prop}

\begin{proof}

(ii) $\Rightarrow$ (i) easily follows from the fact that the natural
norm of
$c_0(\Gamma)$ is metric-$KK^*$.

(i) $\Rightarrow$ (ii) relies on

\begin{Fac} If $X$ has a metric-$KK^*$ norm, there exists a P.R.I.
$(P_\alpha)$ on $X$ such that for any $\alpha < {\rm dens}(X)$, if
$(x,y) \in
X^2$ are such that $P_\alpha(x)=x$ and $P_\alpha(y)=0$, then $\Vert x+y
\Vert
=\max(\Vert x\Vert,\Vert y\Vert)$.
\end{Fac}

Indeed by Lemma 4.2 we know that $X$ is w.c.g. Then Lemma 4.3 shows
that
for all $x \in X$, there is $E_x \subset X^*$ a separable subspace such
that
$\Vert x+y\Vert=\max(\Vert x\Vert,\Vert y\Vert)$ for every $y\in
(E_x)_\perp$. We now use the technique of ([D-G-Z], Lemma VI.2.3):
using the
same notation, we prove along the same lines that if $A \subset X$ and
$B
\subset X^*$ are subsets with density $\leq \aleph$, there exist norm
closed
subspaces $[A]\subset X$ and $[B] \subset X^*$ such that

(i) $A \subset [A]$, $B \subset [B]$.

(ii) ${\rm dens}([A])\leq \aleph$, ${\rm dens}([B]) \leq \aleph$.

(iii) For all $x \in [A]$, $\Vert x\Vert =\sup \{f(x);\ f \in [B],\
\Vert
f\Vert \leq 1\}$.

(iv) For all $x \in [A]$, $E_x \subset [B]$.

(v) For all $f \in [B]$, for all $s \in S$,
$\displaystyle{\sup_{L_s}}\vert
f\vert = \displaystyle{\sup_{L_s \cap [A]}}\vert f\vert$.

\noindent Note that ([D-G-Z], Lemma VI.2.4) shows that $X=[A] \oplus
[B]_\perp$, while the choice of $E_x$ and (iv) shows that $\Vert
x+y\Vert=
\sup(\Vert x\Vert,\Vert y\Vert)$ for all $x \in [A]$ and $y\in
[B]_\perp$. Fact
4.6 now follows by a simple transfinite induction argument, as in
([D-G-Z],
Theorem VI.2.5).

Finally, Proposition 4.5 follows immediately by transfinite induction from
Fact
4.6, since Theorem 2.4 proves it in the separable case and allows us
to
start the induction.
\end{proof}

Theorem 2.4 shows in particular that a separable Banach space has an
equivalent   $LKK^*$ norm if and only if it has an equivalent
metric-$KK^*$ norm, hence the distinction between the two notions is purely
isometric for separable spaces.  Our next two
statements show that it is not so in the non separable case, since certain
spaces are w.c.g. while others are not.

\begin{Thm} Let $K$ be a compact space. The following assertions
are  equivalent:

\noindent
(i) The Cantor derived set of order $\omega_0$ of $K$ is empty.

\noindent
(ii) $C(K)$ is Lipschitz isomorphic to $c_0(\Gamma)$, where $\Gamma$ is
the
density character of $C(K)$.

\noindent
(iii) $C(K)$ admits an equivalent Lipschitz weak-star Kadec-Klee norm.\end{Thm}

\begin{proof} (i) $\Rightarrow$ (ii) was proved in [D-G-Z 2] and the
argument
for the converse can be found in [J-L-S] (Theorem 6.3). The
equivalence between (i) and (iii)  follows easily from the proof of ([La],
Theorem 3.8).
\end{proof}

Our next statement provides the topological condition which allows
``linearizing'' the Lipschitz isomorphism from Theorem 4.7.

\begin{Thm} Let $K$ be a compact space. The following assertions
are  equivalent:

\noindent
(i) $K$ is an Eberlein compact and its Cantor derived set of order
$\omega_0$ is empty.

\noindent
(ii) $C(K)$ is linearly isomorphic to $c_0(\Gamma)$, where $\Gamma$ is
the
density character of $C(K)$.

\noindent
(iii) $C(K)$ admits an equivalent metric weak-star Kadec-Klee norm.\end{Thm}

\begin{proof}

\noindent (i) implies (ii): Since $K$ is Eberlein, $C(K)$ is w.c.g. (see
[D-G-Z], Chapter VI). By compactness,
$K^{(\omega_0)}=\emptyset$ implies that
there is $n$ in $\Ndb$ such that $K^{(n)}=\emptyset$. We proceed by
induction
on $n$. If $n=1$, $K$ is finite and the implication is obvious. Assume
it
holds when $L^{(n)}=\emptyset$ and pick $K$ such that
$K^{(n+1)}=\emptyset$.
We let $L=K'$ and $X=\{ f\in C(K):\ f_{|L}=0 \}$. The space $X$
is clearly isometric to $c_0(K \setminus L)$; while $C(K)/X$ is
isometric to
$C(L)$, and thus isomorphic to a $c_0(\Gamma)$ space by our assumption.
We
observe now that $X$ is complemented in $C(K)$, since any $c_0(I)$ space
is
2-complemented in any w.c.g. space. For checking this, let us call $Y$ a
subspace isometric to $c_0(I)$ of a w.c.g. space $X$. Using the
notation of
([D-G-Z], section  VI.2), we can choose the map $\varphi: X^* \to X^\Ndb$ from
([D-G-Z], Lemma VI.2.3) in such a way that for any $s \in S$ and any $f
\in
X^*$:

\noindent (i) $\displaystyle{\sup_{L_s}}|f|=\sup\{|f(x)|;\ x \in
\varphi(f)
\cap L_s\}.$

\noindent (ii) $\displaystyle{\sup_{Y\cap L_s}}|f|=\sup\{|f(x)|;\ x \in
\varphi(f)
\cap L_s\cap Y\}.$

\noindent (iii) $\overline{{\rm span}\, \varphi(f)\cap Y}=\{x \in
c_0(I);\ {\rm
supp}(x)
\subseteq I_f\}$, where $I_f$ is a countable subset of $I$.

Then ([D-G-Z], Lemma VI.2.4 and Th VI.2.5) provide a P.R.I.
$(P_\alpha)$ on
$X$ such that for all $\alpha \leq {\rm dens}(X)$:

\noindent 1) $P_\alpha(Y)\subseteq Y$.

\noindent 2) There exists $I_\alpha \subseteq I$ such that
$P_\alpha(x)= 1\kern-2.5pt{\rm I}_{I_\alpha}x$ for all $x\in Y$.

\noindent By Sobczyk's theorem, $c_0(\Ndb)$ is 2-complemented in any
separable
super-space. Then we proceed by induction on dens$(X)$: if it is true
for all
w.c.g. $Z$ with dens$(Z)<$ dens$(X)$, we consider $(P_\alpha)$ which
satisfies
1) and 2) above. Since
$(P_{\alpha+1}-P_\alpha)(c_0(I))=c_0(I_{\alpha+1}\setminus I_\alpha)$,
there
is a projection
$$\Pi_\alpha: P_{\alpha+1}(X) \to (P_{\alpha+1}-P_\alpha)(c_0(I))$$
such that $\Vert \Pi_\alpha\Vert \leq 2$. Let
$\Pi'_\alpha=\Pi_\alpha P_{\alpha+1}$ and $\Pi=\sum \Pi'_\alpha$. It is
easily
checked that $\Pi$ is the required projection from $X$ onto $Y$ with
$\Vert
\Pi\Vert \leq 2$.

To conclude the proof of (i) $\Rightarrow$ (ii), we simply observe that
since
$X$ is complemented in $C(K)$, we have that
$$C(K) \backsimeq X \oplus C(L) \backsimeq c_0(K \setminus L)
\oplus c_0(\Gamma).$$

(ii) implies (iii) is clear since the natural norm on $c_0(\Gamma)$ is
metric-$KK^*$.

(iii) implies (i): By Lemma 4.2, any Banach space which has a metric-$KK^*$
norm has a shrinking P.R.I. and thus is w.c.g. The condition
$K^{(\omega_0)}=\emptyset$ follows immediately from Theorem 4.7.
\end{proof}

{\sl Examples 4.9:} There exist ([C-P]; see [D-G-Z], section  VI.8) compact
spaces
such that $K^{(3)}=\emptyset$ (hence $C(K)$ is Lipschitz isomorphic to
$c_0(\Gamma)$) but there is no continuous one-to-one map from
$(B_{C(K)},w)$
to $(B_{c_0(\Gamma)},w )$ and thus no linear continuous injective map
from such
a $C(K)$ to any $c_0(\Gamma)$. Therefore Theorems 2.1 and 2.2 do not extend
to the non separable case. In fact, each compact space
$K$ such that $K^{(\omega_0)}=\emptyset$ but $K$ is not Eberlein provides
an example and some of these are quite simple (see
[D-G-Z], Example VI.8.7).

\section{CHARACTERIZATIONS OF $c_0(\Gamma)$. ADDITIONAL REMARKS}
\medskip
 In this last section we use the above non separable techniques for
characterizing $c_0(\Gamma)$ spaces by showing that they are
the only ${\cal L}^\infty$ spaces which are ``optimally smooth". This leads
in particular to the extension of our main results to
non separable w.c.g. spaces (Corollary 5.2). We also gather some remarks on
the non separable theory. We begin with:

\begin{Prop} Let $X$ be a Banach space such that {\rm
dens}$(X)=\omega_1$.
The following assertions are equivalent:

\noindent (i) $X$ is linearly isomorphic to $c_0(\Gamma)$, with
$|\Gamma|=\omega_1$.

\noindent (ii) $X$ is a ${\cal L}^\infty$ space with an equivalent
metric weak-star Kadec-Klee norm.

\noindent (iii) $X$ is a weakly compactly generated ${\cal L}^\infty$
space
with an equivalent Lipschitz weak-star Kadec-Klee norm. \end{Prop}

\begin{proof} It is obvious that (i) implies (ii) and (iii).

\noindent (iii) implies (i): We use the notation from the proof of
Theorem 4.4.
Through an easy separable exhaustion argument we can ensure that the
spaces
$X_\alpha$ are (uniformly in $\alpha$) ${\cal L}^\infty$ spaces. By
restriction, they have (uniformly in $\alpha$) $LKK^*$ norms.
Hence by Proposition 3.1 they are uniformly isomorphic to $c_0(\Ndb)$.
This
clearly implies (i).

\noindent (ii) implies (iii) follows immediately from Lemma 4.2.
\end{proof}

We can now prove an extension of Theorems 2.1 and 2.2 to certain non
separable spaces. Examples 4.9 above show that it is necessary
to assume that the spaces are w.c.g. On the other hand, the restriction on
the cardinality of $\Gamma$ aims at avoiding
technichalities and it is probably unnecessary.

\begin{Cor} Let $X$ be a weakly compactly generated Banach space, and let
$\Gamma$ be a set with with $|\Gamma|=\omega_1$.
Then:

\noindent (i) If $X$ is Lipschitz isomorphic to a subspace of
$c_0(\Gamma)$, then it is linearly isomorphic to a subspace of
$c_0(\Gamma)$.

\noindent (ii) If $X$ is Lipschitz isomorphic to $c_0(\Gamma)$, then it is
linearly isomorphic to
$c_0(\Gamma)$.
\end{Cor}

\begin{proof} The proof of Theorem 2.1 shows that if $X$ is Lipschitz
isomorphic to a subspace of $c_0(\Gamma)$, then $X$ has an
equivalent $LKK^*$ norm. Indeed the $LKK^*$ property is separably
determined by definition and an easy exhaustion argument shows
that if $E$ is any separable subspace of $X$, there is a separable space
$F$ with $E\subset F\subset X$ and $F$ is Lipschitz
isomorphic to a subspace of $c_0(\Gamma)$. Now (i) follows from Theorem 4.4
and (ii) from Propostion 5.1 and the fact that being a
${\cal L}^\infty$ space is stable under Lipschitz isomorphisms ([H-M]).
\end{proof}

Our next statement provides an extension of ([G-L], Th. IV.1.; see also
[H-W-W] p. 134) to non separable spaces.

\begin{Prop} Let $X$ be a ${\cal L}^\infty$ space with {\rm
dens}$(X)=\omega_1$ which is isomorphic to an $M$-ideal in its bidual
equipped with its bidual norm. Then $X$ is isomorphic to $c_0(\Gamma)$
where
$|\Gamma|=\omega_1$. \end{Prop}

\begin{proof} Since $X$ is an $M$-ideal in $X^{**}$, it is w.c.g. and it
admits a
shrinking
P.R.I. $(P_\alpha)_{\alpha \leq \omega_1}$ by ([F-G], Th. 3). Let
$\lambda \in
\Rdb$ be such that $X$ is ${\cal L}^\infty_\lambda$. For any sequence
$(x_n^*)$ in $X^*$ with $\Vert x_n^* \Vert=1$ and $w^*-\lim x^*_n=0$,
there
exists $\alpha<\omega_1$ such that:

\noindent (a) $P_\alpha^*(x_n^*)=x^*_n$ for every $n \geq 1$.

\noindent (b) $P_\alpha(X)$ is a  ${\cal L}^\infty_\lambda$ space.

\noindent Since $P_\alpha(X)$ is a separable ${\cal L}^\infty_\lambda$
space
which is $M$-ideal in its bidual, we have by ([G-L], Remark 1, p. 261)
that
$d_{BM}(P_\alpha(X),c_0(\Ndb)) \leq M$, where $M=M(\lambda)$ depends
only upon
$\lambda$. It follows that there exists a cluster point to the sequence
$(x_n^*)$ in $(X^{***},w^*)$, say $G$, such that $d(G,X^*) \geq A>0$, where $A$
depends
only on $M$ (that is, on $\lambda$).

If now $(x_n^*) \subset B_{X^*}$ and $w^*-\lim x^*_n=x^*$, with $\Vert
x^*_n-x^*\Vert \geq \varepsilon$, there is, by the above, $G$ in
$B_{X^{***}}$
with $d(G,X^*) \geq A\varepsilon$ and $G=w^*-\displaystyle{\lim_{\cal
U}}(x^*_n-x^*)$ in $(X^{***},w^*)$. Since
$G+x^*=w^*-\displaystyle{\lim_{\cal
U}}x^*_n$, one has $1\geq \Vert G+x^* \Vert=\Vert G\Vert+\Vert x^*\Vert$
and it
follows that $\Vert x^*\Vert \leq 1-A\varepsilon$. Recapitulating, we
have
shown that any separable subspace of $X$ is $A-LKK^*$.
Finally,
Proposition 5.1 yields the conclusion.
\end{proof}

\begin{Rem} 1) It is clear that any quotient space of $c_0(\Gamma)$ has
a
metric-$KK^*$ norm, namely the quotient norm. Therefore Proposition 4.5
shows that Alspach's
theorem [Al]
extends to arbitrary $c_0(\Gamma)$ spaces. That is, any quotient space of
$c_0(\Gamma)$ is isomorphic to a subspace of
$c_0(\Gamma)$, and the isomorphism constant can be made arbitrarily close
to 1. Similarly, Fact 4.6 shows that Johnson-Zippin's theorem [J-Z] extends
to arbitrary $c_0(\Gamma)$ spaces. That is, a ${\cal
L}^\infty$ subspace of a $c_0(\Gamma)$ space is itself isomorphic to a
$c_0(\Gamma_1)$ space.

\medskip
2) Since Lemma 4.3 only uses separable subspaces of $X$, the proofs of
Theorem 4.4 and Proposition 4.5 provide: let $X$ be a Banach space. If for
every
separable subspace $Y$ of $X$, $d_{BM}(Y,\{$subspaces of
$c_0(\Ndb)\})=1,$
then $$d_{BM}(X,\{{\rm subspaces~ of~ c_0(\Gamma)}\})=1.$$

\noindent Now consider a w.c.g. space $X$ with ${\rm
dens}(X)=\omega_1$ and such that every separable subspace of
$X$ is isomorphic to a subspace of
$c_0(\Ndb)$. An argument by contradiction shows the existence
of an upper
bound $M>0$ for the Banach-Mazur distance of any separable subspace of
$X$ to
the subspaces of $c_0(\Ndb)$. Then we get that $X$ is isomorphic to a
subspace
of $c_0(\Gamma)$. Examples 4.9 show that we cannot
dispense with the
assumption ``$X$ w.c.g." in this case.

\medskip

3) An alternative approach to show Proposition 5.3
consists
into proving (with the same notation) that the sequence $(x^*_n-x^*)$
has a
cluster point $G$ in $(X^{***},w^*)$ with $d(G,X^*) \geq A\varepsilon$
for
some constant $A>0$, by extracting first a subsequence which is
$(\varepsilon
/ 2)$-separated, then a further subsequence which is
$(K\varepsilon)$-equivalent to the unit vector basis of $\ell_1$ for
some
constant $K>0$. Indeed, by [L-S], $X^*$ is isomorphic to
$\ell_1(\Gamma)$ and
thus it has the strong Schur property. Now we can pick a $w^*$-cluster
point
$G$ to that subsequence in $(X^{***},w^*)$ to reach our conclusion. The
interest of this alternative route lies in the fact that in the
separable
case, it provides a proof of ([G-L], Th. IV.1) which relies on Proposition 3.1
instead of using Zippin's converse to Sobczyck's theorem ([Z4]).

\medskip
4) It is not difficult to show (using an argument from [A]) that if $X$
has an
equivalent $LKK^*$ norm, then there is an equivalent norm on
$X^{**}$ such that $X$ is an $M$-ideal in $X^{**}$. But this norm is in
general not the bidual norm of its restriction to $X$: indeed it follows
from [La] that for any
$K$
scattered compact set with $K^{(\omega_0)}=\emptyset$, $C(K)$ has an
equivalent $LKK^*$ norm; but such spaces are not in general
w.c.g.
(see Examples 4.9).
\end{Rem}

\bigskip\bigskip\noindent

\centerline{REFERENCES}
\bigskip\noindent
[Ah] I. AHARONI, Every separable metric space is Lipschitz equivalent to a
subset of $c_0$, {\sl Israel J. Math.}, 19 (1974),
284-291.

\medskip\noindent
[Al] D. ALSPACH, Quotients of $c_0$ are almost isometric to subspaces of
$c_0$, {\sl Proc. Amer. Math. Soc.}, 76 (1979), 285-288.

\medskip\noindent
[A-L] D. AMIR, J. LINDENSTRAUSS, The structure of weakly compact sets in
Banach spaces, {\sl Ann. of Math.} 88 (1968), 35-46.

\medskip\noindent
[A] T. ANDO, A theorem on non empty intersection of convex sets and its
application, {\sl J. Approx. Theory}, 13 (1975), 158-166.

\medskip\noindent
[B] S.J. BERNAU, Small projections on $\ell_1(n)$, {\sl Longhorn notes,
The
University of Texas}, 1983/84.

\medskip\noindent
[C-P] K. CIESIELSKI, R. POL, A weakly Lindel\" of function space $C(K)$
without any continuous injection into $c_0(\Gamma)$, {\sl Bull. Pol.
Acad.
Sci. Math.}, 32 (1984), 681-688.

\medskip\noindent
[D-G-Z] R. DEVILLE, G. GODEFROY, V. ZIZLER, Smoothness and renormings in
Banach spaces, {\sl Longman Scientific and Technical}, 1993.

\medskip\noindent
[D-G-Z2] R. DEVILLE, G. GODEFROY, V. ZIZLER, The three space problem for
smooth
partitions of unity and $C(K)$ spaces, {\sl Math. Annal.}, 288 (1990),
613-625.

\medskip\noindent
[Di] J. DIESTEL, Geometry of Banach spaces. Selected topics, {\sl
Springer
Verlag Lecture Notes in Mathematics}, 485 (1975).

\medskip\noindent
[D-G-K] S. DILWORTH, M. GIRARDI, D. KUTZAROVA, Banach spaces which admit a
norm with the uniform Kadec-Klee property,
{\sl Studia Mathematica}, 112,3 (1995), 267-277.

\medskip\noindent
[F] M. FABIAN, Each countably determined Asplund space admits a
Fr\'echet
differentiable norm, {\sl Bull. Austr. Math. Soc.}, 36 (1987), 367-374.

\medskip\noindent
[F-G] M. FABIAN, G. GODEFROY, The dual of every Asplund space admits a
projectional resolution of the identity, {\sl Studia Math.}, 91 (1988),
141-151.

\medskip\noindent
[G-K-L1] G. GODEFROY, N.J. KALTON, G. LANCIEN, L'espace de Banach $c_0$
est
d\'etermin\'e par sa m\'etrique, {\sl Note aux C.R.A.S.}, t. 327,
S\'erie I (1998), 817-822.

\medskip\noindent
[G-K-L2] G. GODEFROY, N.J. KALTON, G. LANCIEN, Szlenk indices and uniform
homeomorphisms, {\sl preprint}.

\medskip\noindent
[G-L] G. GODEFROY, D. LI, Some natural families of $M$-ideals, {\sl Math.
Scand.}, 66 (1990), 249-263.

\medskip\noindent
[G] E. GORELIK, The uniform nonequivalence of $L_p$ and $l_p$, {\sl
Israel J.
Math.}, 87 (1994), 1-8.

\medskip\noindent
[H-W-W] P. HARMAND, D. WERNER, W. WERNER, $M$-ideals in Banach spaces
and
Banach algebras, {\sl Springer Verlag Lecture Notes in Mathematics},
1547
(1993).

\medskip\noindent
[H-M] S. HEINRICH, P. MANKIEWICZ, Applications of ultrapowers to the
uniform
and Lipschitz classification of Banach spaces, {\sl Studia Math.}, 73
(1982),
225-251.

\medskip\noindent
[Jo-Z] K. JOHN, V. ZIZLER, Smoothness and its equivalents in weakly
compactly
generated Banach spaces, {\sl J. of Funct. Ana.}, 15 (1974), 1-11.

\medskip\noindent
[J-L-S] W.B. JOHNSON, J. LINDENSTRAUSS, G. SCHECHTMAN, Banach spaces
determined by their uniform structures, {\sl Geom. Funct. Analysis}, 3
(1996), 430-470.

\medskip\noindent
[J-L-P-S] W.B. JOHNSON, J. LINDENSTRAUSS, D. PREISS, G. SCHECHTMAN, Moduli
of asymptotic convexity and smoothness and affine
approximation, {\sl in preparation}.

\medskip\noindent
[J-R] W.B. JOHNSON, H.P. ROSENTHAL, On weak*-basic sequences and their
applications to the study of Banach spaces, {\sl Israel J.
Math.}, 9 (1972), 77-92.

\medskip\noindent
[J-Z] W.B. JOHNSON, M. ZIPPIN, On subspaces of quotients of $(\sum
G_n)_{l_p}$
and $(\sum G_n)_{c_0}$, {\sl Israel J. Math.}, 13 (1972), 311-316.

\medskip\noindent
[K-O] N.J. KALTON, M.I. OSTROWSKII, Distances between Banach
spaces, {\sl Forum Math.}, 11 (1999), 17-48.

\medskip\noindent
[K-R] N.J. KALTON, J.W. ROBERTS, Uniformly exhaustive submeasures and
nearly additive set functions, {\sl Trans. Amer. Math. Soc.},
278 (1983), 803-816.

\medskip\noindent
[K-W] N.J. KALTON, D. WERNER, Property $(M)$, $M$-ideals and almost
isometric
structure of Banach spaces, {\sl J. reine angew. Math.}, 461 (1995),
137-178.

\medskip\noindent
[K-O-S] H. KNAUST, E. ODELL, T. SCHLUMPRECHT, On asymptotic
structure, the
Szlenk index and UKK properties in  Banach spaces,{\sl Positivity}, 3
(1999), 173-199.

\medskip\noindent
[La]  G. LANCIEN, On uniformly convex and uniformly Kadec-Klee
renormings,
{\sl Serdica Math. J.}, 21 (1995), 1-18.

\medskip\noindent
[L-S] D.R. LEWIS, C. STEGALL, Banach spaces whose duals are isomorphic
to
$l_1(\Gamma)$, {\sl J. Functional Analysis}, 12 (1973), 177-187.

\medskip\noindent
[Li] J. LINDENSTRAUSS, On non linear projections in Banach spaces, {\sl
Michigan J. Math.}, 11 (1964), 268-287.

\medskip\noindent
[L-T] J. LINDENSTRAUSS, L. TZAFRIRI, Classical Banach spaces, {\sl
Springer
Verlag Lecture Notes in Mathematics}, 338 (1973).

\medskip\noindent
[M] V.D. MILMAN, Geometric theory of Banach spaces. II. Geometry of the
unit ball. (Russian) {\sl Uspehi Mat. Nauk} 26 (1971),
6(162), 73-149.
\noindent
English translation:{\sl Russian Math. Surveys} 26 (1971), 6, 79-163.

\medskip\noindent
[P] D. PREISS, Differentiability of Lipschitz functions on Banach spaces,
{\sl J. Funct. Anal.}, 91 (1990), 312-345.

\medskip\noindent
[Z1] M. ZIPPIN, The finite dimensional ${\cal
P}_\lambda$-spaces with small
$\lambda$, {\sl Israel J. Math.}, 39 (1981), 359-364.

\medskip\noindent
[Z2] M. ZIPPIN, Errata to the paper "The finite dimensional
${\cal P}_\lambda$-spaces with small $\lambda$, {\sl Israel J. Math.},
48,2-3
(1984), 255-256.

\medskip\noindent
[Z3] M. ZIPPIN, ${\cal L}^\infty$ subspaces of $c_0$, {\sl unpublished
preprint},
(1973).

\medskip\noindent
[Z4] M. ZIPPIN, The separable extension problem, {\sl Israel J. Math.},
26
(1977), 372-387.

\bigskip\bigskip\noindent
\noindent{\bf Adresses:}

\noindent
G. Godefroy: Equipe d'Analyse, Universit\'e Paris VI, Bo\^ ite 186, 4,
place
Jussieu, F-75252 Paris cedex 05.

\noindent
N.J. Kalton: University of Missouri, Columbia, MO 65211, U.S.A.

\noindent
G. Lancien: Equipe de Math\'ematiques - UMR 6623, Universit\'e de
Franche-Comt\'e, F-25030 Besan\c con cedex.

\end{document}